\documentclass[12pt]{article}
\usepackage{geometry}                
\usepackage{graphicx}
\usepackage{caption}
\usepackage{subcaption}
\usepackage{epstopdf}
\usepackage{listings}
\usepackage{natbib}

\title{La Nova Scientia: Rewriting the History of Operational Research}
\author{Roberto Rossi}

\begin{document}
\maketitle

In 1622, a Dutch army tried to seize the port of Ma\c{c}ao. A Portuguese Jesuit mathematician did the geometry to determine the distance of a Dutch stockpile of gunpowder landed ashore and the angle of elevation at which the cannon should be set. A direct hit turned the tide of the battle and ensured Ma\c{c}ao remained Portuguese. It is hard to not see parallels between this episode mentioned by \cite{citeulike:14287684} and modern Operational Research (OR) case studies.

In his book ``Operational Research in War and Peace: The British Experience from the 1930s to 1970'' Kirby (2003) provides the following definition of OR.

\begin{quote}
Operational Research is the application of the methods of science to complex problems arising in the direction and management of large systems of men, machines, materials, and money in industry, business, and defence.
\end{quote}

\begin{figure}[h]
\centering
\includegraphics[width=8cm]{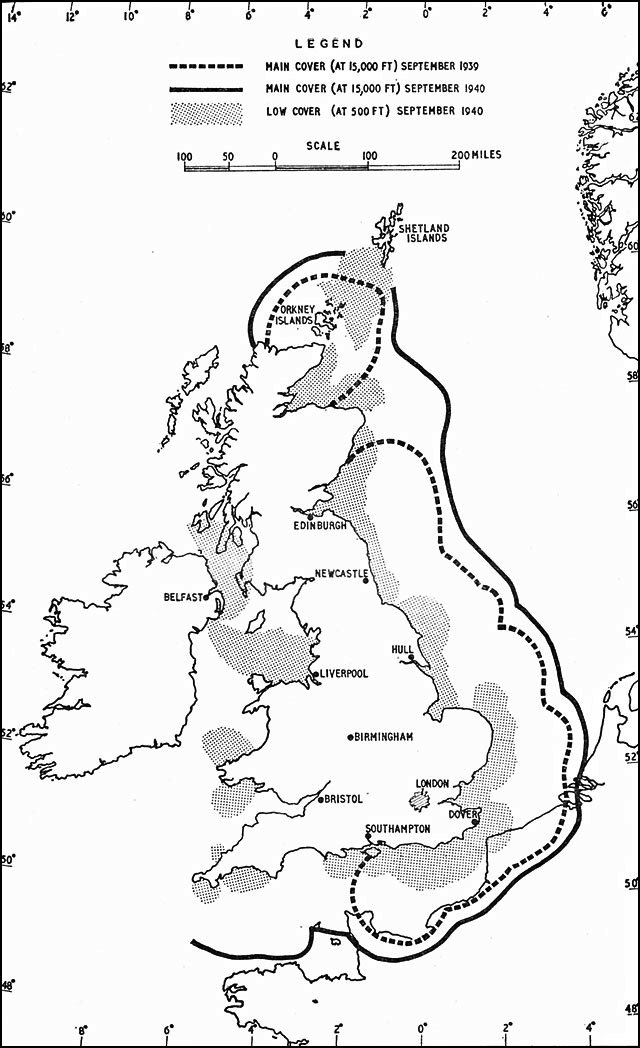}
\caption{Radar coverage along the British coast in the Second World War}
\label{fig:chain_home}
\end{figure}

The  accepted lore is that OR traces its roots back to the First and Second World Wars, when scientific research was used to improve military operations \citep{citeulike:14287763,citeulike:14287762}. Classical examples often discussed include: optimal positioning of radar bases along the UK coastline during the Battle of Britain (Fig. \ref{fig:chain_home}); and military resource allocation problems in the context of the trade-off between aerial antisubmarine-warfare and strategic bombing operations in the Continent \citep{citeulike:14287762}. Of course, cryptanalytic techniques such as Banburismus and Scritchmus developed at Bletchley Park to help break German Kriegsmarine (naval) messages can also be seen as contributions to the field. 

In Chapter 2 of his work \cite{citeulike:14287758} states ``whilst the discipline and practice of OR originated in the later 1930 and came to fruition during the Second World War, the history of warfare is punctuated by attempts to apply some elements of quantitative analysis to understanding the causes of victory and defeat.'' He then goes further to say ``the first overt attempt to apply {\em algebra and numbers} to the battlefield occurred in later sixteenth century.''

As an example of new military practices introduced in the sixteenth century Kirby discusses innovations adopted by Maurice of Nassau, Prince of Orange, in the period between 1585--1625, which included systematising loading and firing of matchlock guns into forty-two sequential moves, so that soldiers could fire more rapidly and in unison; as well as the introduction of the {\em counter march} in which  two ranks of arquebusiers or musketeers alternated their fire to allow reloading. Nevertheless these are qualitative process improvements and Kirby himself defines Maurice as an ``unconscious practitioners'' of OR. 

Kirby's example is only one of many instances of informal applications of OR practices in warfare; another example being the military use of gunpowder in Roberto Valturio's ``De Re Militari'' (1486). A natural question then arises: were OR practices adopted in warfare during the fifteenth and sixteenth centuries always unconscious? In other words, were military practitioners of the sixteenth century simply following common sense, or did some of them follow a principled, scientific approach akin to the one in use today? 

I believe Tartaglia's ``La Nova Scientia'' (1537) represents an exception and should be, in fact, regarded as one of the seminal works in the field of OR, intended as the systematic application of the methods of science to complex problems faced in military operations. 

The fifteenth and sixteenth centuries saw impressive technological advances. Cartography made a quantum leap thanks to the development of perspective. Improved maps, the use of compass and quadrant, as well as the new practice of double-entry bookkeeping allowed merchants to manage ever-expanding trades. Gunpowder and heavy artillery started to be systematically deployed on battlefields and to revolution warfare. 

In his essay ``Metallurgy, Ballistics and Epistemic Instruments'' \cite{citeulike:14287687} investigated how and why a new theoretical ballistic had emerged in the sixteenth centuries. Valleriani observes that that developments in certain practical activities, such as metallurgy, had destroyed the equilibrium between attack and defence strategies that had endured for centuries. As early as 1527 Albrecht D\"urer declared the end of old fortresses, even if they had been readapted. According to his study, the emergence of new theoretical knowledge should be intended as a consequence of an advanced and challenging technological context.

Both \cite{citeulike:14287758} and \cite{citeulike:14287687} date early attempts to apply {\em algebra and numbers} to the battlefield in the late sixteenth century. But how did it all started and what are the essential ingredients to the OR recipe?

\begin{quote}
{\em Painting, cartography and ballistic do not strike us as cutting-edge sciences, but once they were.}
\begin{flushright}
\cite{citeulike:14287687}
\end{flushright}
\end{quote}

The final decline and collapse of the Byzantine empire in the fifteenth century heightened contact between its scholars and those of the west and brought an influx of Neoplatonic scholars. Georgius Gemistus (c. 1355--1452/1454) reintroduced Plato's thoughts to Western Europe during the 1438--1439 Council of Florence, where he met and influenced Cosimo de' Medici to found a new Platonic Academy that focused on the translation of Plato into Latin. As a consequence of the Renaissance revival of Greek mathematics and of the rational tradition of Greek science, we find a proliferation of studies and translations of greek works, including Euclid. 

\begin{figure}[h]
\centering
\includegraphics[width=8cm]{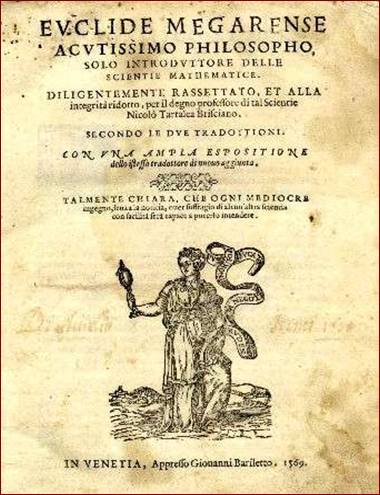}
\caption{Frontespiece of Nicola Tartaglia's traduction of Euclid Elements, 1565}
\label{fig:euclid}
\end{figure}

Euclid's Elements was translated into Arabic in the ninth century. Muslim mathematicians then combined geometry with Hindus arithmetic and algebra and developed new advances. In the twelfth century the work was translated into Latin, making it more accessible to European scholars. Tartaglia knew Euclid, in fact, he went as far as delivering the first translation of Euclid's Elements to Vernacular Italian (Fig. \ref{fig:euclid}).

As discussed by \cite{citeulike:14287684}, Neoplatonism influence was strong at the time of Leon Battista Alberti (1404--1472) and mathematics, as a discipline, tended to be framed in a Neoplatonic context. In his treatise ``On Painting'' Leon Battista Alberti wrote: ``Mathematicians measure the shapes and forms of things in the mind alone and divorced entirely from matter'' \citep{citeulike:14287686}. Reinassance Platonists had been interested in number mysticism rather than real mathematics \citep{citeulike:1853036}. A key ingredient of the OR recipe was missing: {\em the connection between theoretical models and practical applications}. 

An anecdote \citep{citeulike:14287687} discussed by Tartaglia in a letter to Francesco Maria Feltrense della Rovere, Duke of Urbino, suggests that the development of this connection is what motivated the studies discussed in Tartaglia's ``La Nova Scientia.'' In 1532, while he was living in Verona, a friend who was a bombardier asked him: 

\begin{quote}
{\em At which angle the barrel of a cannon should be elevated to achieve the longest possible shot?}
\end{quote}

\begin{figure}[h]
\centering
\includegraphics[width=5cm]{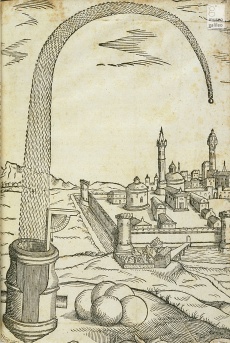}
\caption{Luis Collado. Prattica manuale dell'artiglieria. Milano, 1606, p. 110}
\label{fig:collado}
\end{figure}

Tartaglia did not have expertise in specialised areas connected to military activity. However, having made some calculations, he was able to establish on geometric and algebraic grounds that the maximum range would be achieved if the barrel of a cannon were raised at an angle of 45 degrees above the line of horizon. By answering this question Tartaglia {\em consciously} engaged in an enquiry that today we would have no problem in labelling as OR (Fig. \ref{fig:collado}). 

It is worth remarking that in Tartaglia's time there had been a proliferation of mathematical treatises on topics such as algebra and combinatorics. Both Ramon Lull and Gerolamo Cardano wrote an ``Ars Magna.'' What made Tartaglia's work new (hence the title ``Nova Scientia'') was the application of abstract mathematical models to achieve practical outcomes. This was a revolutionary step that paved the way and inspired - as it is clear from the title of his treatise ``Discorsi intorno a due nuove scienze'' - Galileo's works.

\begin{figure}
\centering
\begin{subfigure}{.45\textwidth}
  \centering
  \includegraphics[width=5cm]{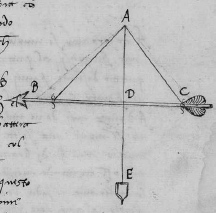}
  \caption{An equilibra, from Leon Battista Alberti's ``Ex Ludis rerum mathematicarum''}
  \label{fig:equilibra}
\end{subfigure}%
\hspace{.01\textwidth}
\begin{subfigure}{.45\textwidth}
  \centering
  \includegraphics[width=5cm]{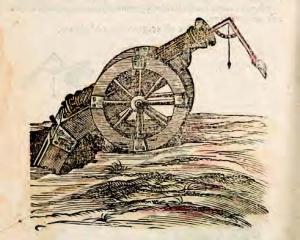}
  \caption{Representation of a cannon positioned at a 45-degree angle of elevation as verified by means of the bombardier's quadrant --- from Tartaglia 1558, epistle, second folio (unnumbered), verso}
  \label{fig:quadrant_gun}
\end{subfigure}
\caption{Equilibra and quadrants}
\label{fig:test}
\end{figure}

However, the connection between theoretical models and applications requires a further ingredient not yet discussed: {\em measurement}. Influenced by Euclid's work, Leon Battista Alberti in his work ``Ex Ludis rerum mathematicarum'' discussed application of trigonometry to surveying \citep{citeulike:14287686}. The methods discussed relied on instruments called equilibra (Fig. \ref{fig:equilibra}). 


Equilibra are a simple extension of the plumb line and could be used to measure angles in everyday activities. By applying basic trigonometrical rules derived from Euclid, these instruments found applications in surveying, e.g. measurements of height of towers and walls. These applications of mathematical tools were different from those devised by Tartaglia: measurement records properties of the World, Tartaglia's applications to ballistics aimed at {\em influencing} the world, not merely observing it.

\begin{figure}[h]
\centering
\includegraphics[width=8cm]{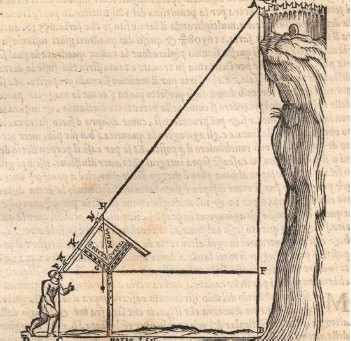}
\caption{Graphic descriptions of the method for measuring heights using the quadrant proposed by Tartaglia. From Tartaglia 1558, third book, 26r}
\label{fig:quadrant_height}
\end{figure}

The first main contribution of Tartaglia consisted in perfecting the quadrant and developing a systematic methodology for its use in ballistic. Tartaglia did not invent the quadrant. A quadrant can be seen as an enhanced equilibra, and equilibra had been in use for a long time. In Tartaglia's work the quadrant was described in two versions developed respectively for aiming cannon (Fig. \ref{fig:quadrant_gun}) --- with a longer leg to be inserted into the gun mouth and a quarter-circle divided into 12 points with plumb line --- and to measure the height and distance of a target. In the latter case, the quadrant has legs of the same length and a shadow square with plumb line (Fig. \ref{fig:quadrant_height}). 

Prior to the era of ballistics the quadrant was chiefly used as an instrument of recording. Before each shot the angle of elevation would be measured and noted. If the shot was successful the annotation would be used to realign the position of the piece of artillery, which would have been lost through recoil \citep[this procedure is discussed in][p. 48-50]{citeulike:14287686}.

Tartaglia's work is not about mere recording and repositioning. He illustrates techniques for estimating distances of target and aiming cannons accordingly. These techniques are based on Euclidean geometry and are illustrated in form of propositions similar to the ones found in Euclid's Elements. Results presented are supported by geometrical (mathematical) reasoning.  

Tartaglia's explicit aim was to create a science that was strictly mathematical and of an axiomatic-deductive nature: he begins with definitions, postulates and axioms, finally propositions and corollaries emerge through a process of deductive reasoning. There is hardly any difference from this structure and that of a modern OR research article.

\begin{figure}[p]
\centering
\includegraphics[width=12cm]{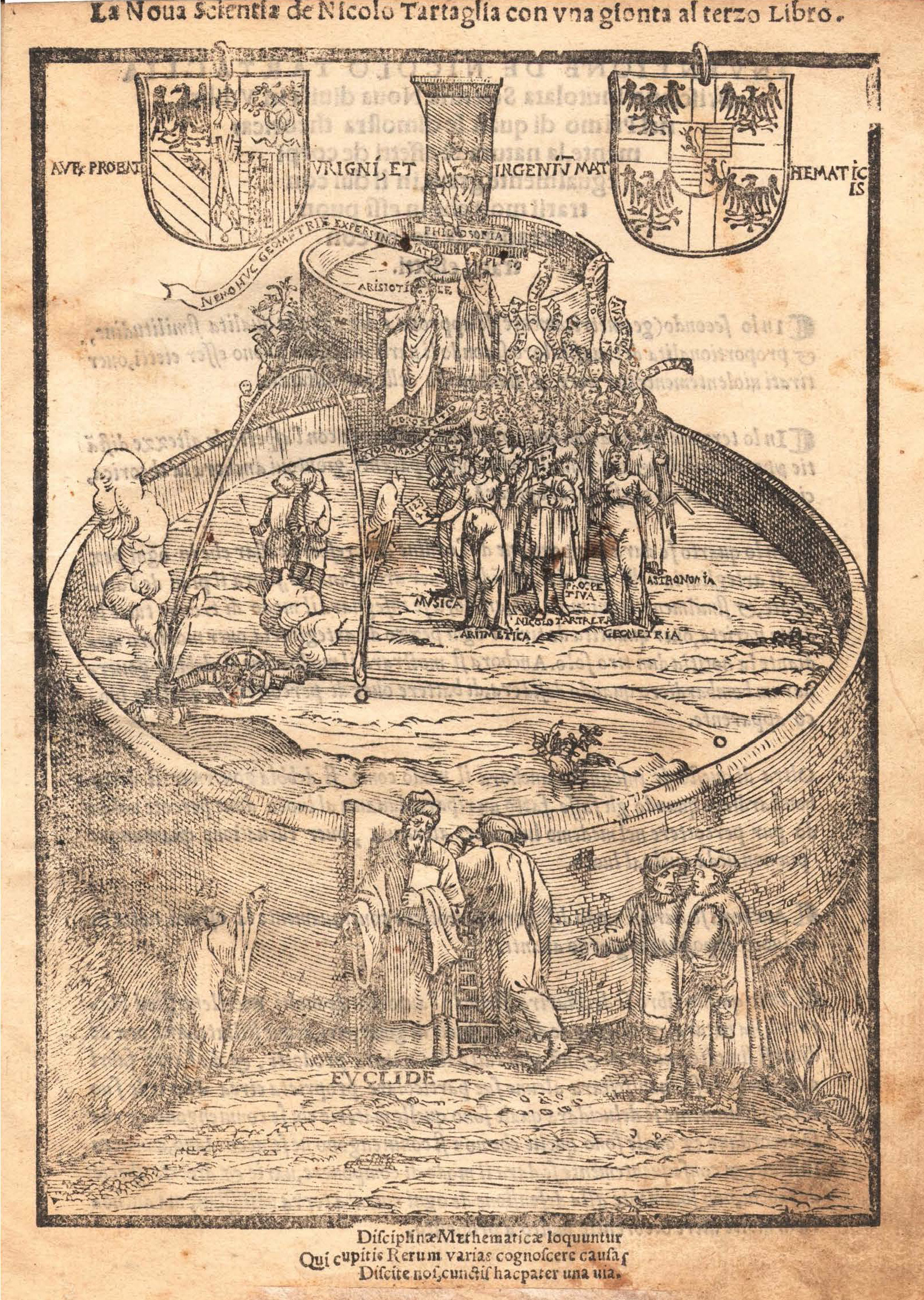}
\caption{Frontispiece of Nicol\`o Tartaglia's Nova scientia (1537). Euclid controls the gate to the castle of knowledge, where a mortar and a cannon are being fired to illustrate path of projectiles. Entry to the inner redoubt requires one to pass through the mathematical sciences, with Tartaglia himself standing among them; within is Philosophy, accompanied by Aristotle and Plato}
\label{fig:frontispiece}
\end{figure}

In the frontispiece of his work (Fig. \ref{fig:frontispiece}) Tartaglia provides a variant of a motto originally found in Luca Pacioli's ``De Divina Proportione.''  In contrast to Pacioli, who refers to the platonic solids, Tartaglia variant of the motto states that the mathematical disciplines (as opposed to the platonic solids) are seen as the only method to understand ``the reasons of things'' and that ``is open to every one;'' this sets a precedent to Galileo's Book of Nature and represents a clear cut from the Hermetic and Neoplatonic traditions, in which knowledge is esoteric.

Tartaglia had to operate within the conceptual framework of Aristotelian physics, as no other viable framework existed. He described the trajectory (``transit'') of projectiles as a sequence of a violent motion followed by a natural motion, which were connected by a circular phase. Tartaglia had limited formal tools to model trajectories mathematically, so he did what every OR person would do: he forced the trajectory to take an approximate form that he could analyse with the formal tools he had. Tartaglia knew that the real trajectory was not made up of  two straight lines joined up by an arc of circumference. He therefore states: ``Nevertheless, that part [of the transit] that is not perceived as being curved is assumed to be straight, and that part that is evidently curved is assumed to be part of the circumference of a circle, as this [assumption] does not influence the argument.'' This approach closely resembles approximation strategies commonly adopted in OR.

\begin{figure}[h]
\centering
\includegraphics[width=8cm]{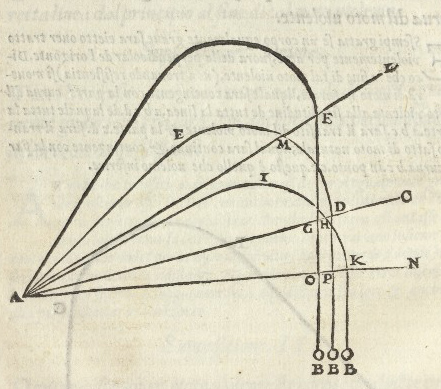}
\caption{Graphic representations of the ``distance of transit'' (distantia del transito) with reference to three different angles of elevation. ``The distance of transit'' is represented by the straight line which joins the initial and final point of the violent motion. From Tartaglia 1558, second book, 11v.}
\label{fig:trajectory}
\end{figure}

On one hand, Tartaglia strived to achieve the greatest possible abstraction of the practical problems he was facing and from the experience of the bombardiers, and thus to construct an exact science based on the Euclidean model. However, Tartaglia was also aware of the applied nature of the {\em Nova Scientia} he was introducing. Some of Tartaglia's arguments - including the one which reveals that the greatest range of shot is achieved when the piece of artillery is raised at an angle of 45 degrees - cannot be explained on a purely geometrical basis, they require {\em observation and experience} in order to be considered valid. This is another aspect that aligns his contribution with existing practices in OR.

Finally, it is worth spending some words on Tartaglia's use of the quadrant as an epistemic instrument \citep{citeulike:14287687}. Valleriani discusses how the annotations made during the fifteenth century by bombardiers represent the beginning of a codified written recording of their experience and practical knowledge. Initially, these annotations could bear meaning only for the bombardier who originally recorded them. The diffusion of the quadrant led to a formalisation of the the measurement process.

\begin{quote}
{\em The annotation regarding the angle of elevation of the piece of artillery is a first step in a process of abstraction, and therefore in theoretical reflection on the bombardier's own actions.}
\begin{flushright}
\citep{citeulike:14287687}
\end{flushright}
\end{quote}

Thanks to the quadrant the bombardier can now describe his activity in a way that is accessible to somebody who is not familiar with his work but possesses the necessary mathematical understanding, such as that of Euclidean geometry, or the necessary physical understanding, such as that of Aristotelian dynamics. The quadrant links theory and practice and enables the transition from bombardier's experience to ballistics as a new theoretical subject. Ultimately, it is this process of {\em abstraction} that makes OR - i.e. the application of the methods of science to complex decision or optimisation problems arising in practical settings - possible.

\bibliographystyle{abbrvnat}
\bibliography{publications}

\end{document}